\documentstyle[12pt]{article}
\newtheorem{defn}{Definition}[section]
\newtheorem{lem}[defn]{Lemma}
\newtheorem{thm}[defn]{Theorem}
\newtheorem{prop}[defn]{Proposition}
\newtheorem{cor}[defn]{Corollary}
\newtheorem{rem}[defn]{Remark}

\newtheorem{assu-def}[defn]{Assumption--Definition}

\newtheorem{nota}[defn]{Notation}
\newtheorem{nota-rem}[defn]{Notation--Remark}

\newcommand{\mod}{\,{\rm mod}\,}
\newcommand{\Z}{{\bf Z}}
\newcommand{\Q}{{\bf Q}}
\newcommand{\R}{{\bf R}}

\newcommand{\Pic}{{\rm Pic}}
\newcommand{\pp}{{\bf P}}

\newcommand{\qed}{$\hfill\,\,\,\diamond$\par\smallskip}
\newcommand{\proof}{{\bf Proof:}\,\,}
\newcommand{\inv}{^{-1}}
\newcommand{\OO}{{\cal O}}
\newcommand{\OS}{\OO_S}

\begin{document}

\title{A connected component of the moduli space of surfaces of general type
with
$p_g=0$}
\author{Margarida Mendes Lopes -- Rita Pardini}
\date{}

\maketitle

\def\theequation{\thesection.\arabic{equation}}
\maketitle

\setcounter{defn}{0} \setcounter{equation}{0}

\begin{abstract}

\noindent Let $S$ be a minimal surface of general type with $p_g(S)=0$ and
such that
the bicanonical map
$\phi:S\to \pp^{K^2_S}$ is a morphism: then the degree of $\phi$ is at most
$4$ by \cite{marg}, and if it is equal to $4$ then $K^2_S\le 6$ by
\cite{untitled}. We
prove that if $K^2_S=6$ and $\deg \phi=4$ then $S$ is a so-called {\em
Burniat
surface} (see
\cite{peters}).  In addition we show that minimal surfaces with $p_g=0$,
$K^2=6$ and
bicanonical map of degree
$4$ form a $4$-dimensional  irreducible connected component  of the moduli
space of
surfaces of general type.
\end{abstract}

\section{Introduction}

Let $S$ be a minimal surface of general type over the complex numbers with
$p_g(S)=0$,
and denote by
$\phi:S\to \pp^{K^2_S}$  the bicanonical map: in \cite{marg}, the first
author has
proven that if $K_S^2\ge 5$, or if $K^2_S= 3, 4$ and
$\phi$ is a morphism, then the degree of $\phi$ is $\le 4$. This result is
made more
precise in \cite{untitled}, where it is proven that if $\deg \phi=4$, then
$K^2_S\le
6$. The latter bound is sharp, as it is shown by the the so-called {\em
Burniat
surfaces} (see
\cite{peters} and \cite{burniat}): these are surfaces of general type with
$p_g=0$, $2\le K^2 \le 6$ whose bicanonical map is $4$-to-$1$ onto a smooth
Del Pezzo
surface. Burniat surfaces  arise as minimal  desingularizations of
$\Z_2\times
\Z_2$-covers of
$\pp^2$ branched on certain arrangements of lines; a direct construction for
the case
$K^2=6$ is given in section \ref{burniat}.

Here we concentrate on the ``borderline case'',  namely $K^2_S=6$. We start
by
showing that {\em all} these surfaces have smooth bicanonical image and
ample
canonical class.  This is an unexpected feature, and may perhaps be related
to the
fact that, although these surfaces have moduli (see theorem
\ref{intromoduli}), they
are expected to be rigid.

The smoothness of the bicanonical image is  the starting point for a very
detailed
analysis of the geometry of these surfaces that enables us to prove the
following:
\begin{thm}\label{introburniat} Let $S$ be a minimal surface of general type
with
$p_g(S)=0$,
$K^2_S=6$ and bicanonical map of degree $4$: then $S$ is a Burniat surface.
\end{thm} This result is also  somehow surprising, since the Burniat
construction is
apparently very  special, and one would not expect it  to include all the
possible
examples.

Theorem \ref{introburniat} also gives us  a good understanding  of the
moduli of the
surfaces we are studying: in fact, using natural deformations of
$\Z_2\times\Z_2$-covers, we are able to prove:
\begin{thm}\label{intromoduli} Minimal surfaces $S$ of general type with
$p_g(S)=0$,
$K^2_S=6$ and bicanonical map of degree $4$ form a $4$-dimensional
irreducible
connected component  of the moduli space of surfaces of general type.
\end{thm}

The plan of the paper is as follows: in section \ref{covers} we explain some
facts on
irregular double covers of surfaces with
$p_g=0$, which are our   main technical tool; in section \ref{burniat}, we
recall the
construction of  Burniat surfaces  and we study their deformations; in
section
\ref{image} we prove that the bicanonical image is smooth: the proof is long
and not
very enlightening, but, as explained above, this is a key result;  in
section
\ref{varie} we collect all the technical facts that we use to prove the main
results
of section \ref{mainres}.
\bigskip

\noindent {\em Acknowledgements:}  Part of this work was carried out
during a stay of the authors at  Warwick University in June 1999:  we are
grateful to the Maths Institute for the wonderful  hospitality and the
stimulating environment we found there; above all we are indebted to Miles
Reid for  many
interesting discussions.   We  are also grateful to Barbara Fantechi for
allowing us to include some unpublished work of hers and the second author
on
deformations of Burniat surfaces and for many useful conversations. Finally,
we
wish to thank  Ciro Ciliberto for  helping us to solve some last minute
doubts.

 This research  has been partly  supported by the italian P.I.N. 1997
 ``Geometria Algebrica, Algebra Commutativa e aspetti computazionali'' and
by Project
PRAXIS XXI 2/2.1/MAT/73/94. The first author is a member of CMAF and of the
Departamento de Matem\'atica da Faculdade de Ci\^encias da Universidade de
Lisboa and
the second author is a member of GNSAGA of CNR.
\bigskip

\noindent{\bf Notations and conventions}: we work over the complex numbers;
all
varieties are assumed to be  compact and algebraic. We do not distinguish
between
line bundles and divisors on a smooth variety,  and we use both the additive
and the
multiplicative notation. Linear equivalence is denoted by $\equiv$. All the
notation
is standard in algebraic geometry; we just recall here the notation for the
invariants of a surface $S$: $K_S$ is the {\em canonical class},
$p_g(S)=h^0(S,K_S)$ is the {\em geometric genus} and $q(S)=h^1(S,\OS)$ is
the {\em
irregularity}.
\section{Irregular double covers and fibrations}\label{covers}
\setcounter{defn}{0} \setcounter{equation}{0}

We describe here the key idea of several   proofs in this paper.

Let $S$ be a smooth surface, let $D\subset S$ be a smooth curve (possibly
empty) and
let $M$ be a line bundle on $S$ such that
$2M\equiv D$. It is well known that there exist a smooth surface
$Y$ and a finite degree $2$ map $\pi:Y\to S$ branched over $D$ and such that
$\pi_*\OO_Y=\OS\oplus M\inv$. We will refer to $S$ as to the {\em double
cover given
by the relation $2M\equiv D$}. The invariants of $Y$ are:
\begin{eqnarray}\label{formule} K_Y^2& =&2(K_S+M)^2\\
\chi(\OO_Y)&= &2\chi(\OS)+\frac{1}{2}M(K_S+M)\nonumber\\
p_g(Y)&=&p_g(S)+h^0(S,K_S+M)\nonumber
\end{eqnarray} If $p_g(S)=q(S)=0$, then the existence of a double cover
$Y\to S$ with
$q(Y)>0$ forces the existence of a fibration
$f:S\to\pp^1$ such that the inverse image via $\pi$ of the general fibre of
$f$ is
disconnected. This is  an old  result of   De Franchis (\cite{defra}), which
is
explained and generalized in several ways in \cite{cetraro}. However, since
these
references are perhaps not widely available, we state it here:
\begin{prop}(De Franchis)\label{defranchis}
 Let $S$ be a smooth surface such that  $p_g(S)=q(S)=0$ and let  $\pi:Y\to
S$ be a
smooth double cover; if
$q(Y)>0$,  then:

i) the Albanese image of $Y$ is a curve $B$;

ii) let $\alpha:Y\to B$ be the Albanese fibration:   there exist a fibration
$g:S\to
P^1$ and a degree $2$ map $p:B\to P^1$ such that\,\,
$p\circ \alpha = g\circ \pi$.
\end{prop}
\proof Denote by $\sigma:Y\to Y$ the involution induced by $\pi$:
$\sigma$ acts on the Albanese variety of $Y$ as multiplication by $-1$,
since
$q(S)=0$. Given   $\eta_1, \eta_2\in H^0(Y,\Omega^1_Y)$,
$\theta=
\eta_1\wedge  \eta_2 $ is a global $2$-form on $Y$ that is invariant for
$\sigma$,
and so it induces an element $\theta'\in H^0(S, K_S)$. Since
$p_g(S)=0$, $\theta'$ vanishes identically, and so does $\theta$. So the
Albanese
image of $Y$ is a curve $B$. The involution $\sigma$ acts on
$Y$ and on
$B$ in a compatible way, and thus the fibration $\alpha:Y\to B$ induces a
fibration
$g:S\to B/\!<\sigma>$. Finally, the quotient curve
$B/\!<\sigma>$ is isomorphic to $\pp^1$, since $q(S)=0$.
\qed

Once constructed such a double cover, sometimes one can  reach a
contradiction either
by showing that the restriction of $M$ to the general fibre of the pencil
$g:S\to\pp^1$ is nontrivial,  and therefore the inverse image via $\pi$ of a
general
fibre of $f$ is connected, or by using the following:
\begin{cor}\label{genus2} Let $S$ be a minimal surface of general type such
that
$p_g(S)=q(S)=0$ and $K_S^2\ge 3$,  and let
$\pi:Y\to S$ be a smooth double cover: then $K_Y^2\ge 16(q(Y)-1)$.
\end{cor}
\proof Since the  statement is of course true for $q(Y)\le 1$, we  assume
that
$q(Y)\ge 2$. By proposition \ref{defranchis}, the Albanese map of $Y$ is a
pencil
$\alpha:Y\to B$ and there exists $g:  S\to\pp^1$ such that $g\circ \pi$ is
composed
with $\alpha$. If
$f$ is the genus of a smooth fibre of
$\alpha$ (and thus of $g$), then by \cite{bea} p. $344$ one has: $K^2_Y\ge
8(q(Y)-1)(f-1)$. If the
  inequality in the statement does not hold, then one has
$f\le 2$.  Since $S$ is of general type, one must have $f=2$. On the other
hand, by
\cite{xiaocan} p. $37$, $S$ has no genus $2$ pencil and so we have a
contradiction.
\qed

Finally, we also exploit this contruction to show the existence of multiple
fibres of
fibrations of $S$, as   explained in the following:
\begin{rem}\label{fibre} Le $S$ be a smooth surface and let $\pi:Y\to S$ be
a smooth
double cover; let $g:S\to\pp^1$ be a fibration such that the general fibre
of $g\circ
\pi$ is not  connected, so that there is  a commutative diagram:

\begin{equation}
\begin{array}{rcccl}
\phantom{1} & Y &\stackrel{\pi}{\rightarrow} & S & \phantom{1} \\
\scriptstyle{g'}\!\!\!\!\!\! & \downarrow & \phantom{1} &
\downarrow & \!\!\!\!\!\!
\scriptstyle{g}
\\
\phantom{1} & B & \stackrel{\bar{\pi}}{\rightarrow} & \pp^1 &
\phantom{1}
\end{array}
\end{equation} where $B$ is a smooth curve of genus $b$ and $\bar{\pi}$ is a
double
cover; if $k$ is the cardinality of the image in $\pp^1$ of the branch locus
of
$\pi$, then $g$ has at least $2b+2-k$ fibres that are divisible by $2$. In
particular, if $\pi$ is unramified, then $g$ has at least $2b+2$ fibres
divisible by
$2$.
\end{rem}
\proof Let
$D\subset
\pp^1$ be the branch locus of $\bar{\pi}$, let $\Delta$ be the branch locus
of $\pi$
and let $D_0=g(\Delta)\subseteq D$; by the commutativity of the above
diagram,
$\pi:Y\to S$ is obtained from
$\bar{\pi}$ by base change and normalization, and thus
$g^*D=\Delta +2 \Delta_0$. Thus
$g^*(D-D_0)\subseteq 2\Delta_0$, i.e. the fibres of $g$ over the points of
$D-D_0$
are divisible by $2$.
\qed

\section{The Burniat construction}\label{burniat}
\setcounter{defn}{0} \setcounter{equation}{0}

We recall briefly the construction of Burniat surfaces with $K^2=6$ (see
\cite{peters} and \cite{burniat}), and we  describe their bicanonical map
and their
small deformations.

Let $\Sigma$ be the blow-up of $\pp^2$ at three distinct non collinear
points
$P_1, P_2, P_3$: we denote by $l$  the pull-back of a line in $\pp^2$,  by
$e_i$  the
exceptional curve corresponding to $P_i$, by $f_i\equiv l-e_i$    the strict
transform of a general line   through $P_i$,
$i=1,2,3$, and   by $e'_i$ the  strict transform of the line joining $P_j$
and $P_k$,
$i\ne j\ne k \ne i$. The subscripts will often be regarded as classes $\mod
3$.  The
$e_i'$'s are disjoint $-1$--curves that also arise  as the exceptional
curves of a
blow-up  map $\Sigma\to\pp^2$: the two blow-ups are related by a quadratic
transformation of $\pp^2$ centered at $P_1, P_2,P_3$. The Picard  group of
$\Sigma$ is the free abelian group generated by the classes of $l, e_1, e_2,
e_3$;
the anticanonical class
$-K_{\Sigma}\equiv 3l-e_1-e_2-e_3\equiv (f_1+f_2+f_3)$ is very ample,  and
the system
$|-K_{\Sigma}|$ embeds
$\Sigma$  as a smooth surface  of degree $6$ in $\pp^6$, the so-called del
Pezzo
surface of degree $6$.

The Burniat surfaces are
$\Z_2\times
\Z_2$-covers of $\Sigma$. Denote by $\gamma_1, \gamma_2,\gamma_3$ the
nonzero
elements of $\Gamma=\Z_2\times\Z_2$ and by $\chi_i\in
\Gamma^*$ the nontrivial character orthogonal to $\gamma_i$; by
\cite{ritaabel}, propositions $2.1$ and $3.1$, to define a smooth
$\Gamma$-cover  $\pi:S\to \Sigma$ one assigns:

 i) smooth divisors
$D_i$,
$i=1,2,3$, such that $D=D_1+D_2+D_3$ is a normal crossing divisor,

ii)  line bundles
$L_1$, $L_2$ satisfying $2L_1\equiv D_2+D_3$, $2L_2\equiv D_1+D_3$.

The branch locus of $\pi$ is $D$: more precisely, $D_i$ is the image of the
divisorial part of the fixed locus of
$\gamma_i$ on
$S$. One has
$\pi_*\OS=\OO_{\Sigma}\oplus L_1\inv\oplus L_2\inv\oplus  L_3\inv$, where
$L_3=L_1+L_2-D_3$ and
$\Gamma$ acts on $L_i\inv$  via the character
$\chi_i$.

To construct a Burniat surface $S$ with $K_S^2=6$, one takes for
$i=1,2,3$  two  smooth   divisors
$m_1^i,m_2^i\in |f_i|$ in such a way that no three of the $m_j^i$'s  have a
common
point and sets:
$D_1=e_1+e'_1+m_1^2+m_2^2$, $D_2=e_2+e'_2+m_1^3+m_2^3$,
$D_3=e_3+e'_3+m_1^1+m_2^1$, $L_1=3l-2e_1-e_3$, $L_2=3l-2e_2-e_1$. By the
above
discussion there exists a smooth $\Gamma$-cover
$\pi:S\to\Sigma$ corresponding to this choice of data, with
$L_3=3l-2e_3-e_2$. The bicanonical divisor
$2K_S=\pi^*(2K_{\Sigma}+D)= \pi^*(-K_{\Sigma})$ is ample, being the
pull-back of an
ample divisor, and thus $S$ is a minimal surface of general type and
$K_S^2=4\frac{1}{4}K_{\Sigma}^2=6$. The invariants of $S$ are:
$\chi(S)=\chi(\pi_*\OS)=1$, $p_g(S)=\sum h^0(\Sigma, K_{\Sigma}+L_i)=0$ and
thus
$q(S)=0$, since $S$ is of general type.

\begin{prop}\label{bicburniat} Let  $S$ be a Burniat surface with $K^2_S=6$:
then
the   bicanonical map of $S$ is the composition of the degree $4$ cover
$\pi:S\to\Sigma$ with the anticanonical  embedding of $\Sigma$ as the smooth
Del
Pezzo surface of degree $6$
 in $\pp^6$.
\end{prop}
\proof Since $p_2(S)=1+K^2_S=7$,  the system
$\pi^*|-K_{\Sigma}|$ is complete,  and the bicanonical map of $S$ is the
composition
of $\pi$ with the anti-canonical embedding of
$\Sigma$ in $\pp^6$. \qed
\bigskip

The last part of this section contains some
  unpublished work  of B. Fantechi and the second author on deformations of
Burniat
surfaces.

\begin{lem}\label{dimdef} Using  the previously introduced notations, one
has:

i)   $h^r(\Sigma, T_{\Sigma}\otimes L_i\inv)=0$, $r\ne 1$, and
$h^1(\Sigma, T_{\Sigma}\otimes L_i\inv)=6$,
$i=1,2,3$;

ii) $h^2(\Sigma, T_{\Sigma}(-\log D_i)\otimes L_i\inv)\le 2$,
$i=1,2,3$.
\end{lem}
\proof When no confusion is likely to arise we omit to write  the space
where
cohomology groups are taken.
 Let $\epsilon:\Sigma\to \pp^2$ be the map that blows down
$e_1,e_2,e_3$; pulling back the Euler sequence on $\pp^2$ one obtains:
\begin{equation}\label{pulltangent} 0\to \OO_{\Sigma}\to
\OO_{\Sigma}(l)^3\to
\epsilon^*T_{\pp^2}\to 0.
\end{equation} As $-L_i+l\equiv -(f_i+e'_{i+1})$, one  has the restriction
sequence:
$$0\to -L_i+l\to \OO_{\Sigma}\to \OO_{f_i}\oplus \OO_{e'_{i+1}}\to 0.$$  The
corresponding long exact sequence gives $h^1(-L_i+l)=1$,
$h^r(-L_i+l)=0$ for $r\ne 1$.
 Tensoring \ref{pulltangent} with $L_i\inv$, we get
$h^r(\epsilon^*T_{\pp^2}\otimes L_i\inv)=3$ if
$r=1$ and zero otherwise, since  $-L_i$ has no cohomology for $i=1,2,3$.  We
have a
short exact sequence:
 \begin{equation}\label{tangente} 0\to T_{\Sigma}\to\epsilon^*T_{\pp^2}\to
\oplus_i
\OO_{e_i}(-e_i)\to 0.
\end{equation} Since   $(-e_i-L_i)e_i=-1$, $(-e_{i+1}-L_i)e_{i+1}=1$ and
$(-e_{i+2}-L_i)e_{i+2}=0$, claim i)  follows by tensoring \ref{tangente}
with
$L_i\inv$ and considering the corresponding long cohomology sequence. Next
we prove
ii) for, say, i=1; by Serre duality,
$H^2(T_{\Sigma}(-\log D_1)\otimes L_1\inv)$ is the dual of
$H^0(\Omega^1_{\Sigma}(\log D_1)(e_2-e_1))\subseteq
H^0(\Omega^1_{\Sigma}(\log
D_1)(e_2))$. For $i=2, 3$, denote by
$\psi_i:\Sigma\to\pp^1$ the morphism given by $|f_i|$ and let
$\psi=\psi_2\times\psi_3:\Sigma \to Q=\pp^1\times\pp^1$; $\psi$ contracts
$e_1$ to a
point $R_1$  and $e'_1$ to a point $R'_1$,  and it is an isomorphism on
$\Sigma-(e_1\cup e'_1)$. Set $M_i=\psi(m^2_i) \in
\OO_Q(1,0)$, $i=1,2$,  and $N=\psi(e_2)\in \OO_Q(0,1)$ and take   $\sigma
\in
H^0(\Omega^1_{\Sigma}(\log D_1)(e_2))$: then $(\psi\inv)^*\sigma$ is a
section of
$\Omega^1_Q(\log\, M_2\!+\!M_1)(N)$ on $Q-\{R_1, R'_1\}$, and thus
$(\psi\inv)^*\sigma\in H^0(Q,\Omega^1_Q(\log\, M_2\!+\!M_1)(N))$, since
$Q$ is nonsingular. The  linear map
$(\psi\inv)^*:H^0(\Omega^1_{\Sigma}(\log
D_1)(e_2))\to H^0(Q,\Omega^1_Q(\log\, M_2\!+\!M_1)(N))$ so defined is
clearly
 injective.  To finish the proof, we remark
$\Omega^1_Q(\log\, M_2\!+\!M_1)(N)\cong \OO_Q(0,1)\oplus \OO_Q(0,-1)$, and
therefore
$h^0(Q,\Omega^1_Q(\log\, M_2\!+\!M_1)(N))=2$.\qed

\begin{prop}\label{locdef} Let $S$ be a Burniat surface with $K^2_S=6$;
then:

i) the Kuranishi family of $S$ is smooth;

ii) every small deformation of $S$ is also a Burniat surface;
\end{prop}
\proof Again, we omit to write the space where cohomology groups are taken,
if no
confusion is likely to arise. We denote by $p:{\cal X}\to B_0$ the family of
natural
deformations of the cover $\pi:S\to
\Sigma$,  defined in section $5$ of
\cite{ritaabel} (for generalizations and a functorial  approach to natural
deformations see also \cite{ritabarb}),  we  let $B\subset B_0$ be the
maximal open
subset over which $p$ is smooth, and we let
$O\in B$ be the point corresponding to $S$.   Notice that for every
$b\in B$
$p\inv(b)$ is a Burniat surface, since
$H^0(D_i-L_j)=H^0( e_i-e_j) =0$ for
$i\ne j$. The base scheme $B$ is smooth and thus, in order to prove i) and
ii) it is
enough to show that the characteristic map
$\rho:T_{B,0}\to H^1(S, T_{S})$ is surjective.
 Given a vector space $V$ with a
$\Gamma$-action, we write   $V^{inv}$  for the  invariant part and
$V^{(i)}$ for the
subspace on which $\Gamma$ acts via the character $\chi_i$; $\Gamma$ acts
both on
$T_{B,0}$ and
$H^1(S, T_{S})$ and  $\rho$ is equivariant with respect to this action. Thus
we have
a decomposition $\rho=\rho^{inv}\oplus
\rho^1\oplus\rho^2\oplus\rho^3$, where $\rho^{inv}:T_{B,O}^{inv}
\to H^1(S, T_{S})^{inv}$ and $\rho^i:T_{B,O}^{(i)}\to H^1(S, T_S)^{inv}$. By
the
definition of natural deformations, we have
$T_{B,0}^{inv}=\oplus_i H^0( D_i)$,
$T_{B,0}^{(i)}=H^0(D_i-L_{i+1})\oplus H^0(D_i-L_{i+2})$ and thus
$T_{B,0}^{(i)}=0$
for $i=1,2,3$ by the remark above. By proposition $4.1$ of \cite{ritaabel},
one has
$H^1(S, T_S)^{inv}=H^1(\Sigma, T_{\Sigma}(-\log D))$ and $H^1(S,
T_S)^{(i)}=H^1(\Sigma, T_{\Sigma}(-\log D_i)\otimes L_i\inv)$,
$i=1,2,3$. So we have to show:

i) $H^1(\Sigma, T_{\Sigma}(-\log D_i)\otimes L_i\inv)=0$ for
$i=1,2,3$;

ii) $\rho^{inv}:\oplus_i H^0(\Sigma, D_i)\to H^1(\Sigma, T_{\Sigma}(-\log
D)$ is
surjective.

\noindent By
\cite{ritaabel}, proposition
$5.2$,
$\rho^{inv}$ is obtained, up to sign, by composing  the restriction map
$\oplus_i
H^0(\Sigma, D_i)\to  \oplus_i H^0(\OO_{D_i}(D_i))$ with  the coboundary map
induced by the sequence:
\begin{equation}\label{logsequence} 0\to T_{\Sigma}(-\log D)\to
T_{\Sigma}\to\oplus_i\OO_{D_i}(D_i)\to 0
\end{equation} Thus ii) follows from the fact that $\Sigma$ is rigid and
$q(\Sigma)=0$.  Replacing
$D$ with $D_i$ in  sequence
\ref{logsequence}, tensoring with $L_i\inv$ and taking cohomology we get the
following  sequence:

$$0\to H^1( T_{\Sigma}(-\log D_i)\otimes L_i\inv)\to H^1(T_{\Sigma}\otimes
L_i\inv)\to $$
\begin{equation}\label{sequence2}
\to  H^1(\OO_{D_i}(D_i-L_i))
\to H^2(T_{\Sigma}(-\log D_i)\otimes L_i\inv)\to 0
\end{equation} Sequence \ref{sequence2} is exact  on the right by lemma
\ref{dimdef}.
 The components of $D_i$ are all smooth rational curves and
$D_i-L_i=3e_i-3e_{i+1}$
has degree
$-3$ on each of them, so that $h^0(\OO_{D_i}(D_i-L_i))=0$ and
$h^1(\OO_{D_i}(D_i-L_i))=8$. Thus
\ref{sequence2} is also exact on the left,  and ii) follows from  lemma
\ref{dimdef},  considering  the dimensions of the vector spaces in sequence
\ref{sequence2}.
\qed

\begin{thm}\label{globdef} Burniat surfaces with $K^2_S=6$ form an
irreducible open
set
 of dimension $4$ of  the moduli space of surfaces of general type.
\end{thm}
\proof As in the proof of proposition \ref{locdef} we consider the family
$p:{\cal
X}\to B$ of smooth  natural deformations of $S$: for every $b\in B$,
$p\inv(b)$ is a Burniat surface and every Burniat surface occurs as a fibre
of $p$.
The image $U$ of $B$ in the moduli space of surfaces of general type is open
by
proposition \ref{locdef}, ii).  Denote by $f:B\to\pp(H^0(\Sigma, D_1))\times
\pp(H^0(\Sigma, D_2))\times \pp(H^0(\Sigma, D_3))$ the restriction to $B$ of
the
projection map;
$f(B)$ is open,
$6$-dimensional,  and the natural map $B\to U$ induces a map
$f(B)\to U$.  Let
$b, b'\in B$ such that there exists an isomorphism $\psi:S\to S'$, where
$S=p\inv(b)$ and
$S'=p\inv(b')$: the  covers
$\pi:S\to\Sigma$ and $\pi':S'\to \Sigma$ are given by the bicanonical map,
and
therefore there exists an automorphism $\bar{\psi}$ of
$\Sigma$ such that $\bar{\psi}\circ \pi=\pi'\circ\psi$.  Conversely, given
$\bar{\psi}\in Aut(\Sigma)$, then the cover
$\pi':S'\to \Sigma$ given by taking base change of
$\pi:S\to \Sigma$ with $\psi$ gives a Burniat surface $S'$ isomorphic to
$S$. So the
fibre of $f(B)\to U$ has a  map with finite fibres onto
$Aut(\Sigma)$ and thus has dimension $2$. As a consequence,
$\dim U=4$. \qed

\section{The bicanonical image}\label{image}
\setcounter{defn}{0} \setcounter{equation}{0}

From now on  we will stick to the following
\begin{assu-def}\label{ipotesi} We denote by $S$ a smooth minimal surface of
general
type with invariants $K^2_S=6$, $p_g(S)=q(S)=0$; we denote by $\phi:S\to
\Sigma=\phi(S)\subset
\pp^6$
 the bicanonical map, which is a morphism by \cite {red}, and assume that
$\deg\phi=4$. The surface $\Sigma$ has degree $6$.
\end{assu-def}
\begin{rem} As explained in section \ref{burniat}, Burniat surfaces with
$K^2=6$
satisfy assumption \ref{ipotesi}.
\end{rem} We  use the notation introduced in section \ref{burniat}.

\begin{thm}\label{smoothimage} Let $\phi:S\to \Sigma$ be  as in
\ref{ipotesi}: then
$\Sigma$ is  the smooth Del Pezzo surface of degree $6$ in  $\pp^6$ (cf.
section
\ref{burniat}).
\end{thm}
\proof The bicanonical image $\Sigma$ is a linearly normal surface of degree
$6$; so, by theorem $8$ of \cite{nagata}, $\Sigma$ is the image of
$\psi:\hat{\pp}\to \pp^6$,   where $\hat{\pp}$ is the blow-up of
$\pp^2$ at points
$P_1, P_2, P_3$ such that $|-K_{\hat{\pp}}|$ has no fixed components,   and
$\psi$ is
given by the system
$|-K_{\hat{\pp}}|$. Thus the $P_i$'s can be   infinitely near, but it is not
possible
that
$2$ of them are distinct and both  infinitely near to the third one. We
denote by
$l$ the pull-back on $\hat{\pp}$ of a general line in $\pp^2$,  by
$e_i$ the exceptional divisor over
$P_i$, and  by $l_i$ a general line through $P_i$, if $P_i$ is not an
infinitely near
point; moreover we write $L$, $L_i$ for the strict transform on $S$ of
$l$, respectively $l_i$.
$\Sigma$ is smooth iff
$P_1, P_2, P_3$ are distinct and not collinear iff $\hat{P}$ contains no
$-2$--curves; in all the other cases,
$\psi$ contracts to rational double points  the $-2$ curves of
$\hat{\pp^2}$, that are
either components  of the $e_i$'s  or  the strict transform of a line
containing all
the $P_i$'s, if  such a line exists. The proof is  a case  by case
discussion of  the
possible configurations of the $P_i$'s that give rise to singular
$\Sigma$'s: in each
case we consider the pull-back of a general hyperplane section of
$\Sigma$ through one of the singular points,  use it to construct an
irregular double
cover
$\pi:Y\to S$ and then obtain a contradiction by means of the techniques of
section
\ref{covers}.
 \bigskip

\noindent{\em Case A:} The points $P_1, P_2, P_3$, not necessarily  all
distinct, lie
on a line
$m$, whose strict transform on $\hat{P}$ is mapped by $\psi$ to a point
$x\in
\Sigma$.

\noindent  Considering the pull-back on $S$ of a hyperplane section of
$\Sigma$
through $x$,  one can write: $2K_S=2L+Z$, where
$Z$ is effective with
$K_SZ=0$. In particular, $h^0(S,L)\ge 3$. Write $Z=2Z'+Z''$, with $Z''$
reduced; we
wish to show $Z''=0$. Since $K_S(K_S-L-Z')$ is even, we have
$(Z'')^2=4(K_S-L-Z')^2\equiv 0\, (\mod 8)$.  So, if $Z''\ne 0$, then, being
reduced,
it contains at least $4$ irreducible  $-2$-curves. On the other hand, $S$
contains at
most
$3$ irreducible
$-2$-curves, since $h^{1,1}(S)=4$ and therefore $Z''=0$. If  $\pi:Y\to S$ is
the
unramified double cover given by the relation
$2(K_S-L-Z')\equiv 0$, then by the double cover formulas \ref{formule} we
have:
$\chi(Y)=2$, $K_Y^2=12$, $p_g(Y)=h^0(S, 2K_S-L-Z')=h^0(S, L+Z')\ge h^0(S,L)=
3$, and
therefore $q(Y)\ge 2$. This contradicts corollary
\ref{genus2}, and thus this case does not occur.
\bigskip

\noindent{\em Case B:} there is no line containing all the $P_i$'s.

\noindent Assume that $P_3$ is infinitely near to $P_2$: there are two
subcases,
according to whether $P_2$ is infinitely near to $P_1$ or not.
\medskip

\noindent{\em Case B1:} $P_2$ is not infinitely near to $P_1$

\noindent We write $2K_S=2L_2+L_1+Z$, where  $Z$ is effective such that
$K_SZ=0$; one
can show that $Z=2Z'$, with $Z'$ effective, by the  argument used  in case
A. Let
$\pi:Y\to S$ be the double cover branched on a general
$L_1$  and given
 by the relation
$2(K_S-L_2-Z')\equiv L_1$; by the double cover formulas \ref{formule}, one
gets
$\chi(Y)=3$,
$p_g(Y)=h^0(S, 2K_S-L_2-Z')=h^0(S, L_1+L_2+Z')\ge 4$ and thus $q(Y)\ge 2$.
By
proposition \ref{defranchis}, the Albanese image of $Y$ is a curve  and
there exists
a pencil $g:S\to\pp^1$ such that $\pi\circ g$ factorizes through  the
Albanese pencil.
Since $\pi$ is branched on $L_1$, $g$ must be the map  given by
$|L_1|$ and thus $g$ has at least $5$ fibres divisible by $2$, by remark
\ref{fibre}. Write $D=\phi^*(\psi(e_1))$ and denote by $\bar{D}$ the strict
transform
of $\psi(e_1)$,  so that $D=\bar{D}+Z$ with $Z$ effective and $K_SZ=0$; one
has
$D^2=-4$, $DK_S= \bar{D}K_S=2$. If $R$ is the ramification divisor of
$\phi$, by
adjunction one has
$K_S=R+\phi^*K_{\Sigma}$, hence $R=3K_S$;  denoting the fibres of
$g$ that are  divisible by $2$ by $2M_i$,
$i=1\dots 5$, we have $R\ge \sum_i M_i$.  Assume that
$\bar{D}$ is reduced,  and thus it  has no common component with $R$: we
have
$2=K_S\bar{D}=\frac{1}{3}R\bar{D}\ge \frac{1}{3}
\bar{D}\sum_i M_i\ge
\frac{5}{6}\bar{D}L_1=\frac{10}{3}$, and thus we have reached a
contradiction. Next
we  assume that $\bar{D}=2E$, with $E$ irreducible such that $K_SE=1$; in
this case
$L_1E=2$ and so, for every $i$,
$M_iE=1$,  and the point $M_i\cap E$ is  smooth for $E$ and  it is  a
ramification
point of the degree $2$ map $\phi|_E:E\to \psi(e_1)$. Thus $p_a(E)\ge 2$ by
the
Hurwitz formula. On the other hand, one has $0=ZD=2EZ+Z^2$ and
$-4=D^2=(2E+Z)^2=4E^2-Z^2$ and thus $E^2\le -1$, $p_a(E)\le 1$. So case B1
does not
occur.
\bigskip

\noindent{\em Case B2:} $P_2$ is infinitely near to $P_1$.

\noindent  As in the previous cases, write $2K_S=3L_1+Z$, where $Z$ is
effective with
$K_SZ=0$. Since $K_SL_1=4$, the index theorem gives either:

 a) $L_1^2=0$ or,

b) $L_1^2=2$.

 In addition
$8=2K_SL_1=3L_1^2+L_1Z$ implies $L_1Z=8$ in case a) and $L_1Z=2$ in case b).
Taking
squares, one gets
$24=4K_S^2=9L_1^2+6L_1Z+Z^2$ and thus $Z^2=-24$ in case a) and
$Z^2=-6$ in case b). The irreducible  components of $Z$ are
$-2$--curves and there are at least two of them, since $-Z^2$ is not  twice
a square.
 On the other hand, notice that the classes
$L$ and
$\phi^*(\psi (e_3))$ span a $2$-dimensional subspace $V$ in
$H^2(\Sigma, \Q)$, since they are both effective and satisfy
$L^2=4$ and $L\phi^*(\psi (e_3))=0$. Recalling that
$h^2(S,\Q)=4$ and that the classes of irreducible  $-2$--curves are
independent and
orthogonal to $V$, one sees that  there are at most $2$ such curves  on $S$.
So, we
denote by $\theta_1$, $\theta_2$ the irreducible $-2$ curves of $S$ and we
write
$Z=a_1\theta_1+a_2\theta_2$, with $a_i>0$. Observe that
$\theta_1\theta_2\ne 0$, since otherwise we would have integral solutions of
$a_1^2+a_2^2=12$, ($=3$ in case b)). Thus
$\theta_1\theta_2=1$, since the intersection form is negative definite on
the span of
$\theta_1$ and $\theta_2$. The equality $Z^2=-24$ ($=-6$ in case b)) can be
rewritten
as
$(a_1-a_2)^2+a_1a_2=12$ ($= 3$ in case b)). If we assume $a_1\ge a_2$, then
the only
solution is  $a_1=4,\, a_2=2$ in case a) and
$a_1=2,\, a_2=1$ in case b). In addition, $L_1\theta_1=2$ in case a),
$L_1\theta_1=1$ in case b) and $L_1\theta_2=0$ in both cases.

Consider now case a) and let $\pi:Y\to S$ be the  double cover  ramified on
a
general  $L_1$ and given by the relation
$2(K_S-L_1-2\theta_1-\theta_2)\equiv L_1$; we have $\chi(Y)=3$,
$p_g(Y)=h^0(S, 2K_S-L_1-2\theta_1-\theta_2)= h^0(S,
2L_1+2\theta_1+\theta_2)\ge 3$
and thus $q(Y)=1$. So we  argue as in case A, and we see that the pencil
$|L_1|$ on
$S$ is induced by the Albanese pencil of $Y$. The curve
$\Delta=\pi^*\theta_1$ is not contained in a fibre of the Albanese pencil of
$Y$ and
it is a smooth rational curve, since
$\theta_1L_1=2$ and $L_1$ is general. Thus we have a contradiction and case
a) is
ruled out.

In case b), we consider the double cover $\pi:Y\to S$ branched on
$L_1+\theta_2$, $L_1$ general, given by the relation
$2(K_S-L_1-\theta_1)\equiv \L_1+\theta_2$; as usual: $\chi(Y)=3$,
$p_g(Y)=h^0(S, 2K_S-L_1-\theta_1)=h^0(S, 2L_1+\theta_1+\theta_2)\ge 3$ and
thus
$q(Y)\ge 1$. As in the previous cases, the Albanese image of $Y$ is a curve
and the
Albanese pencil induces a base point free linear pencil
$|F|$ on
$S$, that satisfies $L_1F=0$; the index theorem applied to $L_1, F$ gives a
contradiction, and the proof is complete.\qed

\begin{prop}\label{ample} The canonical divisor $K_S$ of $S$ is ample and
$\phi$ is
finite.
\end{prop}
\proof By theorem  \ref{smoothimage},  we have $h^2(\Sigma)=h^2(S)=4$ and
thus the
pull-back map $\phi^*:H^2(\Sigma, \Q)\to H^2(S, \Q)$, being injective, is an
isomorphism preserving the intersection form up to multiplication by $4$. If
a curve
$C$ were contracted by
$\phi$, then the class of $C$ in $H^2(S, \Q)$ would be in the kernel of the
intersection form on $H^2(S,\Q)$, contradicting Poincar\'e's duality. So
$\phi$ is
finite and,  as a consequence, $K_S$ is ample.
\qed

\section{Divisors, pencils and torsion of $S$}\label{varie}

\setcounter{defn}{0} \setcounter{equation}{0}

This section collects all the facts needed in the proof of the main theorem
\ref{main}. By theorem \ref{smoothimage}, if
$\phi:S\to\Sigma$ is as in
\ref{ipotesi}, then $\Sigma$ is isomorphic to the blow-up of $\pp^2$ at
three
distinct non collinear points $P_1, P_2,P_3$ and it is embedded in
$\pp^6$ by the anticanonical system. We study in great detail the pull-back
via
$\phi$ of the  exceptional curves and of the  free pencils of $\Sigma$, and
we
produce  a subgroup
$G\simeq
\Z_2^3$ of $Pic(S)$ that plays a very important role in the proof of theorem
\ref{main}.

 Divisors on
$\Sigma$ are denoted as in section
\ref{burniat}.
\begin{lem}\label{pullback} Let $\phi:S\to\Sigma$ be as in \ref{ipotesi} and
let
$C\subset
\Sigma$ be a $-1$--curve: then  either: i) $\phi^*C$ is a smooth rational
curve with
self-intersection $-4$; or ii) $\phi^*C=2E$, where
$E$ is an irreducible curve with $E^2=-1$, $K_SE=1$.
\end{lem}
\proof One has: $(\phi^*E)^2=-4$, $K_S\phi^*E=2$. So, if $\phi^*E$ is
irreducible
then it is smooth rational and we are in case i). Assume that
$\phi^*E$ is reducible: then $\phi^*E=A+B$, with $A$, $B$ irreducible and
such that
$K_SA=K_SB=1$, $K_S$ being ample by proposition \ref{ample}. If $A\ne B$,
then $AB\ge
0$, $A^2+B^2+2AB=-4$ and so, by parity considerations, either one has
$A^2=B^2=-3$,
$AB=1$ or, say,
$A^2=-3$, $B^2=-1$, $AB=0$. In both cases, the matrix
$\left(\begin{array}{cc} A^2 & AB \\ AB & B^2
\end{array}\right)$ is negative definite, and thus the classes of $A$ and
$B$ span a
$2$-dimensional subspace $V$ of $H^2(S,\Q)$. The projection formula:
$C\phi_*D=D\phi^*C$, for  $C$ and $D$  curves on $\Sigma$ and $S$
respectively,
implies that
$V$ and
$\phi^*(<E>^{\perp})$ are orthogonal subspaces. By Poincar\'e's duality,
$H^2(\Sigma,
\Q)=<E>\oplus^{\perp}<E>^{\perp}$ and thus
$H^2(S,\Q)=\phi^*<E>\oplus^{\perp}\phi^*(<E>^{\perp})$, since, as we have
already
remarked in the proof of proposition \ref{ample}, $\phi^*$ is an isomorphism
preserving the intersection form up to multiplication by $4$. Thus
$V\subseteq
\phi^*<E>$,  contradicting the fact that $V$ has dimension $2$. So we must
have
$A=B$ and we are in case ii).
\qed

\begin{lem}\label{miyaoka} If $S$ is a surface as in assumption
\ref{ipotesi}, then
$S$ does not contain $2$ smooth disjoint rational curves with
self-intersection
$-4$.
\end{lem}
\proof Let $r$ be the cardinality of a set  of smooth  disjoint rational
curves
$D$  on
$S$ such that $D^2=-4$; by \cite{miyaoka}, 2.1, one has the following
inequality:
$r\frac{25}{12}\le c_2(S)-\frac{1}{3}K^2_S=4$, namely $r\le 1$.
\qed
\begin{prop}\label{elliptic} Let  $\phi:S\to\Sigma$ be  as \ref{ipotesi} and
let
$e_i, e'_i\subset \Sigma $, $i=1,2,3$, be  defined as  in section
\ref{burniat}: then
for $i=1,2,3$ there exist irreducible curves
$E_i, E'_i\subset S$ such that $\phi^*e_i=2E_i$, $\phi^*e_i'=2E_i'$ and
$E_i^2=(E'_i)^2=-1$,
$K_SE_i=K_SE'_i=1$.
\end{prop}
\proof By lemmas \ref{pullback} and \ref{miyaoka}, we may assume that there
exist
irreducible curves $E_2, E_3, E'_1, E'_3$ on $S$ such that
$E_i^2=(E'_i)^2=-1$, $K_SE_i=K_SE'_i=1$ and
$\phi^*e_2=2E_2$, $\phi^*e_3=2E_3$, $\phi^*e'_1=2E'_1$,
$\phi^*e'_3=2E'_3$ and that $\phi^*e_1$, $\phi^*e'_2$  either are of the
same type
or they are smooth rational curves. So assume that $\phi^*e_1=R$ is a smooth
rational
curve. Writing
$F_i=\phi^*f_i$ for $i=1,2,3$, one has:
 $2K_S\equiv F_1+F_2+F_3\equiv F_1+R+2E'_3+2E'_1+2E_2\equiv R+2F_1+2E'_1
$. Let
$\pi:Y\to S$ be the double cover corresponding to the  relation
$2(K_S-F_1-E'_1)\equiv R$:
$Y$ is a smooth surface with invariants
$\chi(Y)=2$, $K^2_Y=14$,
$p_g(Y)=h^0(S, 2K_S-F_1-E'_1)= 3$ (see formulas \ref{formule}). The last
equality
follows from the fact that $\phi$ maps
$F_1$ to a conic and
$E'_1$ to a line intersecting the conic transversely  at one point.
Therefore we
have $q(Y)=2$  and the result follows from remark \ref{genus2}. The proof
for
$E'_2$ is similar.
\qed

\begin{nota}\label{eta}  Let $\phi:S\to \Sigma$ be as in
\ref{ipotesi}. By theorem \ref{smoothimage}, $\Sigma$ is the blow-up of
$\pp^2$ at
three non colllinear points and we use the notation of section \ref{burniat}
for
divisors on $\Sigma$; in addition, we  write
$F_i=\phi^*f_i$ and we denote by $g_i:S\to \pp^1$ the morphism given by
$|F_i|$, for
$i=1,2,3$. Often we  use residue classes ($mod\, 3$) for the subscripts. For
instance, the pencil
$g_i$ has two reducible double fibers, that we write as
$2E_{i+1}+2E'_{i+2}$ and $2E_{i+2}+2E'_{i+1}$. We set:
$\eta_i=E_{i+1}+E'_{i+2}-E_{i+2}-E'_{i+1}$, $i=1, 2, 3$, and
$\eta=K_S-(\sum_j E_j+\sum E'_j)$.
\end{nota}

\begin{prop} Let $\phi:S\to \Sigma$ be as in  \ref{ipotesi};  let
$\eta, \eta_1, \eta_2, \eta_3\in \Pic(S)$  be defined as in \ref{eta} and
let $G$ be
the subgroup of $\Pic(S)$ generated by these elements. Then $G=\{0,  \eta_1,
\eta_2,
\eta_3, \eta, \eta+\eta_1, \eta+\eta_2,
\eta+\eta_3\}$, $\eta_1+\eta_2+\eta_3=0$, and $G\simeq\Z_2^3$.
\end{prop}
\proof  It is obvious from the definitions that $2\eta=2\eta_i=0$ and
$\eta_1+\eta_2+\eta_3=0$. In addition,  $\eta=K_S-\sum_j(E_j+E'_j)\ne 0$ and
$\eta+\eta_i= K_S-(E_i+E'_i+2E'_{i+1}+2E_{i+2})\ne 0$,  because
$p_g(S)=0$. Finally, $\eta_i\ne 0$, $i=1,2,3$ by \cite{bpv} lemma (8.3),
chap. III.
So $G$ consists  precisely of the $8$ elements listed above.
\qed

\begin{lem}\label{dim} If $S$ is as in assumption \ref{ipotesi}, then:

i)
$h^0(S, K_S+\eta)=h^0(S, K_S+\eta_i)=1$, $h^1(S, K_S+\eta)=h^1(S,
K_S+\eta_i)=0$,
$i=1,2,3$;

ii) $h^0(S, K_S+\eta+\eta_i)=2$, $h^1(S, K_S+\eta+\eta_i)=1$,
$i=1,2,3$;

iii) if $\tau \in Pic(S)$ is such that $2\tau=0$ and $h^0(S, K_S+\tau)\ge
2$, then
$\tau =\eta+\eta_i$ for some $1\le i \le 3$.
\end{lem}
\proof First of all we remark that if $\tau\in Pic(S)$ satisfies $2\tau=0$,
$\tau\ne 0$, then
$1=\chi(K_S+\tau)=h^0(K_S+\tau)-h^1(K_S+\tau)$,  and therefore
$K_S+\tau$ is effective. Now let
$\tau\in Pic(S)$ be such that $2\tau=0$  and
$h^0(S, K_S+\tau)\ge 2$, and   write $|K_S+\tau|= Z+|M|$, where $Z$ and
$|M|$ are the fixed  and the moving part, respectively. The curves
$2Z+2M$ belong to the bicanonical system
$|2K_S|=\phi^*|-K_{\Sigma}|$,  and thus  $|M|=\phi^*|N|$, where $|N|$ is a
linear
system of $\Sigma$ without fixed components
 such that
$-K_{\Sigma}-2N$ is effective. The only possibility is  $|N|=|f_i|$  for
some
$i=1,2,3$. In turn, this  corresponds to $\tau=\eta+\eta_i$, since
$K_S+\eta+\eta_i=F_i+E_i+E'_i$ and $h^0(S, 2(E_i+E'_i))=1$. In particular,
$h^0(S, K_S+\eta+\eta_i)= 2$.\qed

\begin{lem}\label{fasci} If $S$ is a surface as in assumption \ref{ipotesi},
then
$|F_1|$,
$|F_2|$ and $|F_3|$ are the only irreducible base point free pencils of
$S$.
\end{lem}
\proof Let $D$ be the cohomology  class of a base point free pencil of $S$:
then $D$
lies in the nef cone  $NE(S)\subset H^2(S,\R)$ and satisfies $D^2=0$.
Conversely,
given  $D\in  NE(S)$ with $D^2=0$ there is at most one irreducible pencil of
$S$
whose class is proportional to $D$.

 As we have seen in the proof of corollary.
\ref{ample},
$\phi^*: H^2(\Sigma,
\Q)\to H^2(S,
\Q)$ is an isomorphism preserving the intersection form up to multiplication
by $4$;
in addition,  integral classes both on $S$ and on
$\Sigma$ are algebraic because $p_g(S)=p_g(\Sigma)=0$, and therefore
$NE(S)=\phi^*NE(\Sigma)$.  Now,   $NE(\Sigma)$ is spanned by the classes of
$f_1, f_2, f_3, l, l'$, where $l'$ is the pull-back of a conic in
$\pp^2$ through the fundamental points $P_1, P_2, P_3$,  and so $D$ is equal
to the
class of $f_1$,
$f_2$ or $f_3$.
\qed

\begin{lem}\label{doublefibres}  Let $S$ be as in  \ref{ipotesi} and let
$g_i:S\to
\pp^1$ be as in notation \ref{eta}, $i=1,2,3$; then:

i) the multiple fibres of $g_i$ are double fibers and their number is $\ge
2$ and
$\le 4$;

ii) if $g_i$ has $4$ double fibres, then $E_i$ and $E'_i$ are smooth
elliptic curves.
\end{lem}
\proof We recall that $g_i$ has at least $2$ double fibres, namely
$2E_{i+1}+2E'_{i+2}$ and $2E'_{i+1}+2E_{i+2}$, (see proposition
\ref{elliptic} and notation \ref{eta}).  Let
$mD\in |F_i|$, with $m>1$; since $E_iF_i=E'_iF_i=2$, one has
$m=2$ and   $D$ intersects both $E_i$ and $E'_i$ transversely at smooth
points. The
irreducible curves $E_i$ and
$E'_i$ of arithmetic genus $1$   are mapped by
$\phi$
$2$--to--$1$ onto the smooth rational curves $e_i$ and $e'_i$, and the maps
$E_i\to e_i$ and $E'_i\to e'_i$ are ramified at  the  point  $DE_i$,
respectively
$DE'_i$. So, by the Hurwitz formula, there are at most
$4$ double fibres, and in that case $E_i$ and $E'_i$ are smooth.
\qed

\begin{prop}\label{Fcanonico} Let $S$ be a surface as in assumption
\ref{ipotesi}
and  let $F_i\in |F_i|$ be a general curve, $i=1,2, 3$; if $i\ne j$, then
$F_j|_{F_i}=K_{F_i}$.
\end{prop}
\proof We  show that $F_3|_{F_1}= F_2|_{F_1}=K_{F_1}$. Notice  that
$2K_S=F_1+F_2+F_3=F_1+2E'_3+2E_1+2E'_2+2E_1$, and consider the double cover
$\pi:Y\to
S$  branched on a smooth $F_1$ and given by the relation
$2(K_S-2E_1-E'_3-E'_1)\equiv F_1$; by the formulas \ref{formule}, the
invariants of
$Y$ are $\chi(Y)=3$, $K^2_Y=20$, $p_g(Y)=h^0(S, 2K_S-2E_1-E'_3-E'_2)$. To
give a
lower bound  for $p_g(Y)$ we observe that
$|2K_S-2E_1-E'_3-E'_2|=|(F_1+2E_1)+E'_2+E'_3|=|\phi^*l+E'_2+E'_3|\supseteq
\phi^*|l|+E'_2+E'_3$ (see section \ref{burniat} for the  notation) and
thus
$p_g(Y)=h^0(S, 2K_S-2E_1-E'_3-E'_2)\ge 3$ and $q(Y)\ge 1$. By  proposition
\ref{defranchis}, the Albanese pencil on $Y$ is the pull-back of a pencil
$|F|$ on
$S$ such that $\pi^*F$ is disconnected for $F$ general. Since $\pi$ is
branched on a
curve of $|F_1|$, it follows that
$FF_1=0$ and therefore $|F|=|F_1|$. In addition, if $F_1$ is general  then
$\pi^*F_1$ is the unramified double cover of $F_1$ given by the relation
$2(K_S-2E_1-E'_3-E'_2)|_{F_1}\equiv 0$; since $\pi^*F_1$ is disconnected,
the line
bundle
$(K_S-2E_1-E'_3-E'_2)|_{F_1}=(K_S-2E_1)|_{F_1}=(K_S-F_3)|_{F_1}=
(K_S-F_2)|_{F_1}$ is
trivial.\qed

\begin{prop}\label{etares} Let $S$ be as in  \ref{ipotesi}; for $i=1, 2, 3$
let
$F_i\in |F_i|$ be a  general curve and let
$G_i=\{\tau\in G: \tau|_{F_i}=0, \}$: then
$G_i=\{\eta_i, \eta+\eta_{i+1},\eta+\eta_{i+2}\}$.
\end{prop}
\proof We prove the lemma for $G_1$. One has $\eta_1\in G_1$ by  definition.
Moreover, using lemma \ref{Fcanonico}, it is easy to show that
$\eta|_{F_1}=\eta_2|_{F_1}=\eta_3|_{F_1}=(E_1-E'_1)|_{F_1}$, so we only need
to show
$\eta|_{F_1}\ne 0$. Notice that $K_S+F_1+\eta+\eta_1=
2F_1+E'_1+E_1=2K_S-E_1-E'_1$.
Therefore
$H^0(S,K_S+F_1+\eta+\eta_1)$ is isomorphic to the kernel of the restriction
map
$H^0(S,2K_S)\to H^0(E_1+E'_1, 2K_S|_{E_1+E'_1})$. Since
$|2K_S|$ embeds $E_1+E'_1$ as a pair of skew lines, it follows that
$h^0(S,K_S+F_1+\eta+\eta_1)=3$. Next we restrict
$K_S+F_1+\eta+\eta_1$ to $F_1$ and get:
$0\to H^0(S,K_S+\eta+\eta_1)\to H^0(S,K_S+F_1+\eta+\eta_1)\to
H^0(F_1,K_{F_1}(\eta))\to  H^1(S,K_S+\eta+\eta_1)$. Using lemma
\ref{dim}, it follows that $h^0(F_1,K_{F_1}(\eta))\le 2$ and so
$\eta|_{F_1}$ is nontrivial.
\qed

\section{The main results}\label{mainres}
\setcounter{defn}{0} \setcounter{equation}{0}

This section is devoted to proving of the following:
\begin{thm}\label{main} Let $S$ be a smooth minimal surface of general type
with
invariants
$p_g(S)=q(S)=0$, $K^2_S=6$; if the bicanonical map $\phi:S\to
\Sigma\subset \pp^6$ has degree $4$, then $S$ is a Burniat surface.
\end{thm}
\noindent and
\begin{thm}\label{main2} Smooth minimal surfaces of general type $S$ with
$K^2_S=6$,
$p_g(S)=0$ and bicanonical map of degree $4$ form a  $4$-dimensional
irreducible
connected component of the moduli space of surfaces of general type.
\end{thm}
\noindent {\bf Proof of theorem \ref{main2}:} Let  $\cal M$ be the the
moduli space of
 surfaces  of general type with
$p_g=0$ and $K^2=6$, and let
$\cal{Y}\subset\cal M$  be the subset of surfaces such that the bicanonical
map
has degree $4$: by theorem
\ref{main} and propositions
\ref{bicburniat} and
\ref{locdef}, $\cal{Y}$ is  open in $\cal M$. In addition (cf. \cite{marg}),
$\cal Y$
coincides with  the subset of $\cal M$ consisting of surfaces such that the
bicanonical map has degree $\ge 4$. In order to
show that
$\cal Y$ is also closed,    it is enough to prove the following:
 let  $B$ be an irreducible  curve and let  $f:{\cal X}\to B$ be a smooth
family, such that for every $t$ the fibre $X_t$ is a minimal surface of
general
type and  the bicanonical map $\phi_t:X_t\to\pp^N$ is a generically finite
morphism; then there exists $m$ such that $\deg\phi_t\ge m$ for every $t\in
B$,  with equality holding except  for finitely many points $t\in B$. Up to
normalizing $B$ and restricting to an open subset, we may assume that there
exists $\Phi:{\cal X}\to B\times \pp^N$ such that
$\Phi|_{X_t}=\phi_t$ for every $t\in B$. Denote by $\cal Y$ the image of
$\cal X$
with the reduced  scheme structure: the restriction of the projection
${\cal Y}\to B$
is a flat morphism, since $B$ is smooth of dimension $1$ and $\cal Y$ is
irreducible.
It follows that the fibres
$Y_t$ have constant degree $d$ in $\pp^N$. For every $t\in B$, let $Y'_t$
the reduced
scheme structure underlying $Y_t$: then  one has
$4K^2=\deg\phi_t\deg Y'_t$, and thus $\deg\phi_t\ge m=\frac{4K^2}{d}$, with
equality
holding iff $Y_t$ is generically reduced.\qed
\bigskip

\noindent {\bf Proof of theorem \ref{main}:} Since the  proof is  long, we
break  it
into four steps. We      use the notations introduced in sections
\ref{burniat} and
\ref{varie}. In addition, we denote by
$\pi_i:Y_i\to S$ the unramified double cover given by
$\eta+\eta_i$, for $i =1,2,3$. By the formulas \ref{formule} and lemma
\ref{dim}, we have $p_g(Y_i)=2$, $q(Y_i)=1$; we denote by
$\alpha_i:Y_i\to B_i$ the Albanese pencil.

 {\bf Step 1:} {\em Up to a permutation of $\{1,2,3\}$,  the pencil
$g_{i-1}\circ \pi_i:Y_i\to \pp^1$ is composed with
$\alpha_i:Y_i\to B_i$.}

\noindent By proposition
\ref{defranchis}, the Albanese pencil $\alpha_i:Y_i\to B_i$ arises in the
Stein
factorization of $g\circ \pi_i$ for some base point free pencil
$g:S\to\pp^1$.  By lemma \ref{fasci}, there is
$s_i\in\{1,2,3\}$  such that $g=g_{s_i}$. Notice that $s_i\ne i$, since  by
proposition
\ref{etares} the general curve of
$\pi_i^*|F_j|$ is connected if and only if $i=j$. To prove the claim, we
have to show
that $i\mapsto s_i$ is a permutation of $\{1,2,3\}$. Assume by contradiction
that,
say, $s_2=s_3=1$ and denote by $p: Z\to S$ the unramified
$\Z_2\times\Z_2$-cover  with  data $L_1=\eta_1$, $L_2=\eta+\eta_2$,
$L_3=\eta+\eta_3$
 (see section \ref{burniat},  or \cite{ritaabel} proposition 2.1).  One has
$q(Z)=\sum_i h^1(S, L_i\inv)=2$ by lemma \ref{dim}; we denote by
$\alpha:Z\to A$ the
Albanese map.   If
$\sigma_i$ is the element of $\Z_2\times \Z_2$ that acts trivially on
$L_i\inv$, then, for $i=2,3$, the surface $Z/\!<\sigma_i>$ can be naturally
identified
with
$Y_i$;
 we denote by
$p_i:Z\to Y_i$ the projection map and by $p_{i*}:A\to B_i$ the homomorphism
induced
by $p_i$. Notice that
$p_{2*}\times p_{3*}:A\to B_2\times B_3$ is an isogeny, since $H^1(Z,
\OO_Z)\simeq H^1(S,
\eta+\eta_2)\oplus H^1(S,\eta+\eta_3)\simeq p_2^*H^1(Y_2,
\OO_{Y_2})\oplus p_3^*H^1(Y_3,
\OO_{Y_3})$. Since the pencil $g_1\circ p$ is composed with  both
$p_{2*}\circ
\alpha$ and $p_{3*}\circ
\alpha$, the Albanese image of $Z$ is a curve $B$ of genus $2$ and
$g_1\circ
p=\bar{p} \circ\alpha$, where $\bar{p}:B\to \pp^1$ is a $\Z_2\times
\Z_2$-cover. By
the Hurwitz formula,
$\bar{p}$ is   branched exactly over $5$ points of $\pp^1$, since in a
$\Z_2\times
\Z_2$--cover of smooth curves the inverse image of a branch point consists
of $2$
simple ramification points. Arguing as in the proof of remark \ref{fibre},
one sees
that the fibres of
$g_1$ over the branch points of
$\bar{p}$ are double, but this contradicts lemma
\ref{doublefibres}.

{\bf Step 2:} {\em The general  $F_i$ is hyperelliptic for
$i=1,2,3$.}

\noindent We show that the general $F_1$ is hyperelliptic. We have seen that
the
pencil
$g_1\circ\pi$ is composed with the Albanese map
$\alpha_2:Y_2\to B_2$ and that  $g_3\circ\pi_2$  also has disconnected
fibres. The
Stein factorization of $g_3\circ\pi_2$ is
$Y_2\stackrel{g}\to C \stackrel{\psi}\to \pp^1$ where $g$ has  connected
fibres, $C$
is a smooth curve and $\deg\psi=2$.
 Notice that
$C\cong \pp^1$, since $q(Y_2)=1$ and $g$ is not the Albanese pencil.  Denote
by
$\tilde{F_1}$    a general fibre of $\alpha$ and by $\tilde{F_3}$ a general
fibre of
$g$. From $F_1F_3=4$ it follows that $\tilde{F_1}\tilde{F_3}=2$. So the
linear
system $|\tilde{F_3}|$ cuts out a $g^1_2$ on the general
$\tilde{F_1}$, and thus the general $F_1$ is hyperelliptic.

{\bf Step 3:} {\em The Galois group $\Gamma$  of $\phi:S\to\Sigma$ is
$\Z_2\times\Z_2$}

\noindent For $i=1,2,3$,  denote by $\gamma_i$ the involution on $S$ that
induces the
hyperelliptic involution on the general $F_i$; the $\gamma_i$'s are regular
maps,
since $S$ is minimal, and they belong to $\Gamma$ by proposition
\ref{Fcanonico}.
Consider the involution
$\tilde{\gamma_1}:Y_2\to Y_2$ inducing the hyperelliptic involution on the
general
$\tilde{F_1}$: by construction  $\tilde{\gamma_1}$ maps each $\tilde{F_3}$
to itself,
and the restriction of $\alpha$ to
$\tilde{F_3}$ identifies $\tilde{F_3}/\!<\tilde{\gamma_1}>$  with $B_2$.
 Since  $\pi_2|_{\tilde{F_i}}:\tilde{F_i}\to\pi_2(\tilde{F_i})\in|F_i|$ is
an
isomorphism compatible with the action of $\tilde{\gamma_1}$ and
$\gamma_1$ for $i=1,3$, this implies that $\gamma_1\ne \gamma_3$. In the
same  way
one shows
$\gamma_i\ne \gamma_j$ for $i\ne j$ and thus $\Gamma=\{1,\gamma_1,
\gamma_2, \gamma_3\}$.

{\bf Step 4:} {\em $S$ is a Burniat surface}

\noindent By step 1,  for each $i=1,2,3$ the map $g_{i}\circ \pi_{i+1}$ is
composed
with the Albanese pencil $\alpha_{i+1}:Y_{i+1}\to B_{i+1}$ and  thus, by
remark
\ref{fibre} and lemma \ref{doublefibres}, $g_i$ has precisely $4$ double
fibres. The
double fibres are  $2(E_{i+1}+E'_{i+2})$,
$2(E'_{i+1}+E_{i+2})$, and
$2M^i_1=\phi^*m^i_1$, $2M^i_2=\phi^*m^i_2$, where $m^i_1, m^i_2 \in |f_i|$.
If we
denote by $D$ the total branch locus of $\phi$, then
$D\supseteq D_0=\sum_i(e_i+e'_i+m^i_1+m^i_2)$. By \cite{ritaabel}
proposition
$3.1$,
$D$  is a normal crossing divisor, since $S$ is smooth, and therefore no
three of the
$m_j^i$ have a common point. Applying the Hurwitz formula to a general
bicanonical
curve yields: $-K_{\Sigma}D=18=-K_{\Sigma}D_0$ and thus $D=D_0$, since
$-K_{\Sigma}$ is ample. As in section \ref{burniat}, we denote by $D_i$ the
image of
the divisorial part of the fix locus of
$\gamma_i$, so that $D=D_1+D_2+D_3$. By \cite{ritaabel} proposition $3.1$,
$D_i$ is smooth for every
$i=1,2,3$, so  there is a permutation $i\mapsto s_i$ of $\{1,2,3\}$ such
that
$D_i\supset m_1^{s_i}+m_2^{s_i}$; in addition, the quotient of a general
$F_i$ by
$\gamma_i$ is rational and therefore $D_if_i=4$.  One concludes that for
$i=1,2,3$
$D_i= e_i+e'_i+m_1^{s_i}+m_2^{s_i}$ and $s_i\ne i$. Finally, the quotient of
a general
$F_{i+2}$ by $\gamma_i$ is the elliptic curve $B_{i+1}$ (cf. step 3) and
thus
$D_if_{i+2}=2$.  So one gets  $s_i=i+1$ and
$S$ is obtained precisely as  explained in section
\ref{burniat}.\qed

\bigskip

\begin{tabbing} 1749-016 Lisboa, PORTUGALxxxxxxxxx\= 56127 Pisa, ITALY \kill
Margarida Mendes Lopes               \> Rita Pardini\\ CMAF \> Dipartimento
di
Matematica\\
 Universidade de Lisboa \> Universit\a`a di Pisa \\ Av. Prof. Gama Pinto, 2
\> Via
Buonarroti 2\\ 1649-003 Lisboa, PORTUGAL \> 56127 Pisa, ITALY\\
mmlopes@lmc.fc.ul.pt
\> pardini@dm.unipi.it
\end{tabbing}

\end{document}